\documentclass[12pt,reqno,a4paper]{amsart}
\usepackage{amsmath,amssymb,amsfonts,amscd}
\usepackage[mathscr]{eucal}
\usepackage{hyperref,xcolor} 
\usepackage{comment,textcomp}

\newtheorem{theorem}{Theorem}
\newtheorem{lemma}{Lemma}
\newtheorem{proposition}{Proposition}

\theoremstyle{definition}
\newtheorem{definition}{Definition}
\newtheorem{example}{Example}
\newtheorem{remark}{Remark}

\title[Non-linear homomorphisms...]
{Non-linear homomorphisms of algebras of functions are induced
by thick morphisms}
\author{Hovhannes~M. Khudaverdian}
\address{Department of Mathematics,  
University of Manchester, Manchester,   UK}
\email{khudian@manchester.ac.uk}
\address{%
}

\begin{document}
\begin{abstract}

 In 2014, Voronov introduced the notion 
of thick morphisms of (super)manifolds 
as a tool for constructing $L_{\infty}$-morphisms 
of homotopy Poisson algebras. 
Thick morphisms generalise ordinary 
smooth maps, but are not maps themselves. 
Nevertheless, they induce pull-backs 
on $C^{\infty}$ functions. 
These pull-backs are in 
general non-linear maps 
between the algebras 
of functions which are 
so-called ``non-linear homomorphisms''.
By definition, this means that their differentials 
are algebra homomorphisms in the usual sense.
The following conjecture was formulated:
an arbitrary non-linear
homomorphism of algebras of
smooth functions is generated
by some thick morphism.
We prove here this conjecture
in the class of formal functionals. 
In this way, we extend the 
well-known result for smooth 
maps of manifolds and algebra 
homomorphisms of $C^{\infty}$ 
functions and, more generally, 
provide an analog of classical
``functional-algebraic duality''
in the non-linear setting.
\end{abstract}

\maketitle


\def\A {{\bf A}} 
\def\AA {{\mathfrak A}} 
\def\B {{\cal B}}
\def\C {{\bf C}}
\def\CC {{\cal C}}
\def\Cl {{\tt \hbox{Cliff}}}
\def\E {{\bf E}}
\def\EE {{\cal E}}
\def\F {{\cal F}}
\def\FF {{\cal F}}
\def\G {\Gamma}
\def\GG {{\cal G}}
\def\H {{\bf H}}
\def\K {{\bf K}}
\def\L {{\mathfrak L}}
\def\M {{\cal M}}
\def\N {{\bf N}}
\def\R {{\bf R}}
\def\bS {{\bf S}}
\def\SS {{\mathfrak S}}
\def\Tr {{\rm Tr\,}}
\def\V {{\cal V}}
\def\X {{\bf X}}
\def\XX {{\cal X}}
\def\Y {{\bf Y}}
\def\Z {{\bf Z}}

\def\a {\alpha}
\def\ac {{\bf a}}
\def\b {{\bf b}}
\def\bs {{\bf s}}
\def\c {{\bf c}}
\def\d {\delta}
\def\dist {{\tt \hbox{distance}}}
\def\e {{\bf e}}
\def\f {{\bf f}}
\def\finish {\blacksquare}
\def\g {{\bf g}}
\def\grad {{\rm grad\,}}
\def\h {\hbar}
\def\k {{\bf k}}
\def\l {{\bf l}}
\def\m {{\bf m}}
\def\n {{\bf n}}
\def\p {\partial}
\def\pb {{\bf p}}
\def\pt {{\bf pt}}
\def\q {{\bf q}}
\def\r {{\bf r}}
\def\s {\sigma}
\def\t {{\bf t}}
\def\tS {{\tilde \Sigma}}
\def\td {\tilde}
\def\v {{\bf v}}
\def\vare {\varepsilon}
\def\x {{\bf x}}
\def\y {{\bf y}}
\def\w {\omega}
\def\finish {${\,\,\vrule height1mm depth2mm width 8pt}$}

\section {Introduction}

 A map $\varphi\colon\,\, M\to N$ defines the linear
  map 
         \begin{equation}\label{first}
		 \varphi^*\colon\,\, 
		 C^{\infty}(N)\to C^{\infty}(M)\,,
	 \end{equation}
which is homomorphism of algebras of functions.  
In 2014 Ted Voronov introduced the  notion of a {\it thick morphism}
(see \cite{thickfirst}, \cite{thicktheor}) of manifolds, which
generalises  ordinary maps.  A thick morphism
 defines  a non-linear map 
    $\Phi^*\colon\,\, C^{\infty}(N)\to C^{\infty}(M)$. 
  This notion  
 provides a natural way to construct
 $L_\infty$ morphisms for homotopy Poisson algebras
(see \cite{thickfirst},\cite{thicktheor}, 
and \cite{thickkoszul} and also Appendix A).
The notion of thick morphisms turns out to be also 
related with  
quantum mechanis and the construction of
spinor representation (see \cite{tvoscil} and
\cite{thickaction}).
   The pull-back $\Phi^* \colon\,\, 
		 C^{\infty}(N)\to C^{\infty}(M)$ 
corresponding to a thick morphism is not in general
     a homomorphism of algebras 
    (just because it is non-linear). 
     However as it 
 was proved by Voronov,
the differential of this non-linear map is a  usual pull-back.
This motivated him to define so called {\it non-linear homomorphisms}.
\def\B {{\bf B}}
\begin{definition}(Th.Voronov, see \cite{thicktheor})
	\label{defofnonlinearhom1}
	Let $\A,\B$ be two algebras.
	A map $L$ from an algebra $\A$ to an algebra $\B$
	is  called {\it a non-linear homomorphism} 
	if 
	at an arbitrary element of algebra $\A$
	its derivative  is a homomorphism
	of the algebra $\A$  to the algebra $\B$.
\end{definition}
One can say 
that a thick morphism induces a non-linear homomorphism
of algebras of functions in the same way as a usual morphism
$\varphi$ induces usual (linear) homomorphism \eqref{first}.
 A natural  question was 
formulated in \cite{thicktheor}:  
is it true that every non-linear algebra 
 homomorphism between algebras of
 smooth functions arises from a 
 thick morphism as the pull--back?
 Note that the pull-backs by thick morphisms 
 are formal mappings of algebras. 
 Hence in the above definition of
 non-linear homomorphisms one can consider formal maps only.
 We prove here this conjecture for formal maps
(''formal functionals'').

 The structure of the paper is as follows.
We recall the construction of thick morphisms,
and we define a class of {\it formal functionals}
which are induced by thick morphisms.  
We recall the 
 proof of  Voronov's result that
the functional induced  by a thick morphism
 is a non-linear homomorphism 
(see \cite{thickfirst} and  \cite{thicktheor} for detail).
Then we show that  the converse  implication also holds.
In Appendix A we briefly  discuss the relation of thick morphisms
with $L_\infty$ morphisms of homotopy Poisson algebras.
In appendix B we recall some useful polarisation formulae.

\smallskip

{\bf Acknowledgment.} I am grateful to Th.Th.Voronov not 
only for continuous help during the work on the paper, 
but also for possibility to learn
thick morphisms first-hand.
I am deeply grateful to A.S.Schwarz. He expressed firm  
 belief  that thick morphisms can be formulated in terms of
  functionals which are non-linear homomorphisms.
   This encouraged me to prove Theorem \ref{2},
   the main result of this paper.  
  I am also grateful to A.Verbovetsky for
   many useful comments.  
 
     Part of this work was done during my visit to Lyon in autumn
  2019. I thank O.Kravchenko, C.Roger and Th.Stroble for hospitality.
 This work was partially supported by 
 the LABEX MILYON (ANR-10-LABX-0070) of
 Universit\'e de Lyon, within the program "Investissements d'Avenir"
 (ANR-11-IDEX- 0007) operated by the French National Research Agency (ANR).

\section {Thick morphisms and non-linear functionals}

      Consider two manifolds $M$ and $N$.
  %
   We denote by $x^i$ 
   local coordinates on  $M$,
   and we denote $y^a$ 
   local coordinates on  $N$.
      To define  the thick morphism
      $\Phi\colon\,\,M\Rrightarrow N$
we consider   
a function, $S=S(x,q)$,  
where $x$ is the point on $M$ and $q$ is covector
in $T^*N$.  We supposee that $S=S(x,q)$ is a formal
function, power series over $q$:
               $$
	   S=S(x,q)=S_0(x)+S_1^a(x)q_a+
	   S_2^{ab}(x)q_bq_a+
	   S_3^{cba}(x)q_cq_bq_a+\ldots=
	       $$
   \begin{equation}\label{action1}
	  S_0(x)+S_1^a(x)q_a+S_+(x,q)
	   \,, {\rm where}\,\,
	   S_+(x,q)=
	   \sum_{k=2}^\infty S^{a_1\dots a_k}q_{a_1}\dots q_{a_k}\,,	   
   \end{equation}
coefficients $S_{k}^{a_1\dots a_k}(x)$ are usual smooth functions on $x$.

A formal function $S(x,q)$ is called {\it generating function
 of thick morphism}.  

 \begin{remark}
 In fact $S(x,q)$ is geometrical object which transforms non-
 trivially under changing of local coordinates 
(see for detail \cite{thicktheor}). Here and below we consider
only local coordinates $x^i$ on $M$ and 
$y^a$ on $N$. 
\end{remark}

To generating function  $S(x,q)$ corresponds thick morphism
$\Phi=\Phi_S\, \colon\,\,M\Rrightarrow N$
which is defined in the following way:
it defines pull-back $\Phi_S^*$ such that 
to every smooth function $g(y)\in C^\infty(M)$
corresponds a function

\begin{equation}\label{defofthickmorphism1}
 f(x) =\Phi^*(g)=g(y)+S(x,q)-y^aq_a\,
     \end{equation}
where $y^a=y^a(x),q_b=q_b(x)$ are chosen in a way that
\begin{equation}\label{defofthickmorphism2}
y^a={\p S(x,q)\over \p q_a}\,,\quad
q_b={\p g(y)\over \p y^b}\,.
\end{equation}

\begin{remark}
	Conditions \eqref{defofthickmorphism2} 
	imply  that left hand side of equation 
	\eqref{defofthickmorphism1}
  does not depend on $y^a$ and $q_b$:
	\begin{equation*}
	{\p \over \p y^a}
		\left(
 g(y)+S(x,q)-y^aq_a
		\right)=0\,,\,
	{\p \over \p q_b}
		\left(
 g(y)+S(x,q)-y^aq_a
		\right)=0\,.
	\end{equation*}	
\end{remark}

In the special case if $S(x,q)=S^a(x)q_a$ 
$y^a=S^a(x)$  and $\Phi^*g$ is the usual pull-back corresponding to the
map $y^a=S^a(x)$:
	 \begin{equation}\label{usualpullback1}
 f(x) =\Phi^*(g)=g(y)+S(x,q)-y^aq_a=g\left(S^a(x)\right)\,,
	 \end{equation}
and this pull-back corresponds to the usual morphism $y^a=S^a(x)$.

In the general case  (if action $S(x,q)$ is not linear over $q$) 
maps 
\eqref{defofthickmorphism1} and \eqref{defofthickmorphism2} become
formal maps. They become  formal power series in $g$ 
(see for details also equation \eqref{defofthickmorphism2a} below).  
Namely equation \eqref{usualpullback1} defines the formal
functional $L(x,g)$ on  
$C^\infty(N)$  such that  
		\begin{equation}\label{formalfunctional1}
			L(x,g)= 
			L_0(x,g)+
                          L_1(x,g)+
			L_2(x,g)+\dots=
			\sum L_k(x,g)\,,\, 
	\left(g\in C^\infty(N)\right)
      \end{equation}
where  every summand $L_k(x,g)$  takes values in smooth functions on $M$
and it has an order  $k$ in $g$: 
$L_k(x,\lambda g)=\lambda^kL_k(x,g)$.  We suppose that
	    \begin{equation}\label{generalisedfunctions}
		    L_k(g)=\int L(x,y_1,\dots,y_k)
	    g(y_1)\dots g(y_k)
		    dy_1\dots dy_k\,, 
	    \end{equation}
the kernel  $L(x,y_1,\dots,y_n)$ of the functional  $L_k(x,g)$ 
can be  generalised functions. 

\begin{definition}\label{defofformalfunctionals}
We denote by $\A$ the space of all formal functionals 
	which have appearance \eqref{formalfunctional1}.	
We denote by $\A_k$ the subspace of functionals which have
   order $k$ on $g$, ($k=0,1,2,\dots $).
\end{definition}	   

       For arbitrary functional 
	 $L(x,g)\in \A$ 
	 (see equation \eqref{formalfunctional1})
functional $L_k(x,g)$ 
is projection of functional $L(x,g)$ on subspace $\A_k$.   
We sometimes
denote this  projection by $[L(x,g)]_k$
	 \begin{equation}\label{formalfunctional1a}
			L(x,g)= 
			\sum L_k(x,g)\,,\quad 
		L_k(x,g)=[L(x,g)]_k\,.
	 \end{equation}
It is useful to denote  by $\A_{\geq k}$ ($\A_{\leq k}$) 
the subspace of functionals
which have order bigger  or equal than  $k$ 
(less or equal than $k$),
\begin{equation}\label{lessthank}
       \A_{\geq k}=\oplus_{i\geq k}A_a\,,\quad		    
   \A_{\leq k}=\oplus_{0\leq i\leq k}A_i\,,
\end{equation}	
and we say that two functionals  $L_1, L_2\in\A$ 
 coincide up to the order $k$ if 
$L_1-L_2\in \A_{\geq k+1}$.     
We will write in this case that
		\begin{equation*}\label{equalityoffunctionals}
		    L_1(g)=L_2(g) ({\rm mod}\, \A_{k+1})	   
		     \end{equation*}
Explain how 
every formal generating function $S(x,q)$,
 (see equation 
\eqref {action1}) defines 
   thick morphism
 $\Phi_S$, i.e. how $S(x,q)$  defines a
 map $\Phi_S^*(g)$  which is 
 a  formal functional in $\A$.
  Functional 
	    $\Phi_{S(x,q)}^*(g)$ defines
    non-linear pull-back, assigning to every
   smooth function $g\in C^\infty(N)$ a formal sum 
   of smooth functions 
   $\left[\Phi_{S(x,q)}^*(g)\right]_k$, ($k=0,1,2,\dots$).
   \begin{equation}\label{defofthickmorphism1a}
	    \Phi_{S(x,q)}^*(g)=
	    \sum
	    \left[\Phi_{S(x,q)}^*(g)\right]_k=
	    \left[\Phi_{S(x,q)}^*(g)\right]_0+
	    \left[\Phi_{S(x,q)}^*(g)\right]_1+\dots\,,
   \end{equation}
 where $\left[\Phi_{S(x,q)}^*(g)\right]_k$ is component of
 the functional
	    $\Phi_{S(x,q)}^*(g)$ which  has order 
	    $k$ in $g$ 
	(see equation \eqref{formalfunctional1a}).
We will explain how to calculate this map recurrently 
step by step and we will write explicitly the
results of   calculations of its first components. 
(See Propositions \ref{mapforthickmorphisms} and
\ref{expressionforthickmorphism})

As it was mentioned above a map $y^a=y^a(x)$ in 
	equation \eqref{defofthickmorphism2} 
	has to be viewed as a formal sum of smooth maps
	depending on $g$:
	    \begin{equation}\label{defofthickmorphism2a}
		    y^a(x)=y^a(x,g)=\sum y_k(x,g)=  
		     y^a_0(x)+y^a_1(x,g)+\dots=
	    \end{equation}
Here every term $y^a_k(x)=y^a_k(x,g)$
is a smooth map  of order $k$ in $g$:
            $$
	    y^a_k(x,\lambda g)=\lambda^k
(x,g)\,.
	    $$
We will show how to calculate  
map \eqref{defofthickmorphism2a} step by step recurrently, 
and we will write the expressions for calculating
first few components of this formal map 
(see Proposition \ref{mapforthickmorphisms} below).

One can see from equations \eqref{action1}
and \eqref{defofthickmorphism2}
that initial term $y^a_0(x)$ in equation
\eqref{defofthickmorphism2a} is equal to
	    \begin{equation}\label{formalmapforinitialterm}
	    y_0^a(x)=\left[
		  {\p S(x,q)\over \p q_a}
		    \right]_{q=0}=S^a_1(x)\,,
	    \end{equation}
and every next term $y^a_{k+1}(x)=y^a_{k+1}(x,g)$ in 
\eqref{defofthickmorphism2a} is expressed recurrently via previous
terms $\{y^a_0(x),\dots, y^a_k(x)\}$:
	    \begin{equation}\label{formalmapforarbitraryterm}
	y^a_{k+1}=    
	\left[
{\p S(x,q)\over \p q_a}\big\vert_
{q_a={\p g(y)\over \p y^a}\big\vert_
	{y^a=y^a_{\leq k}(x)}}
\right]_{k+1}\,.
           \end{equation}
Here $y^a_{\leq k}(x)=\sum_{i\leq k} y_i(x)$ 
according to equation \eqref{lessthank}, and
 $[\,\,]_{r}$ means $r$-th component of the map
 (see expansion  \eqref{formalfunctional1a}).

 We have already expression \eqref{formalmapforinitialterm}
 for initial component $y_0(x)$
of map $y^a(x)$ in equation \eqref{defofthickmorphism2a}.
Write down expression for  next components
$y_1^a(x)$ and $y^a_2(x)$ of this map. We have
$$
y^a_{1}=    
	\left[
{\p S(x,q)\over \p q_a}\big\vert_
{q_a={\p g(y)\over \p y^a}\big\vert_{y^a=y^a_{\leq 0}(x)}}
\right]_{1}=
	$$
	\begin{equation}\label{formalmapforfirstterm}
\left[
	\left(
	S^a(x)+2S^{ab}(x){\p g(y)\over \p y^a}
	\right)
_{y^a=S^a(x)}
\right]_1=2S^{ab}(x)g_b^*(x) \,, 
	\end{equation}
       and
	$$
	y^a_{2}=    
	\left[
{\p S(x,q)\over \p q_a}\big\vert_
{q_a={\p g(y)\over \p y^a}\big\vert_{y^a=y^a_{\leq 1}(x)}}
\right]_{2}=	
    $$
     $$
\left[
\left(
    S^a(x)+
   \sum_{j\geq 1}(j+1)S^{a{b_1\dots b_j}}(x)
   {\p g(y)\over \p y^{b_1}}\dots
   {\p g(y)\over \p y^{b_j}}
   \right)
   \big\vert_{y^a=y^a_0(x)+y^a_1(x)}
\right]_{2}=
   $$
  $$
    = \left[
\left(
    S^a(x)+
   2S^{ab}(x)
   {\p g(y)\over \p y^{b}}
        +
   3S^{abc}(x)
   {\p g(y)\over \p y^{b}}
   {\p g(y)\over \p y^{c}}
   \right)
   \big\vert_{y^a=S^a_0(x)+2S^{ab}g_b^*(x)}
\right]_{2}=
$$
   \begin{equation}\label{formalmapforsecondterm}
	   =3S^{abc}(x)g^*_b(x)g^*_c(x)
	   +4S^{ab}(x)S^{cd}(x)g^*_{bc}(x)g_d^*(x)
   \end{equation}
where  in equations 
\eqref{formalmapforfirstterm} and   
\eqref{formalmapforsecondterm}  we used notations 
\begin{equation}\label{usefulnotation1}
	   g^*(x)=g\left(y^a\right)
	    \big\vert_{y^a=S_1^a(x)}\,,\quad
	   g_a^*(x)={\p g(y)\over \p y^a}
	   \big\vert_{y^a=S^a_1(x)}\,,\quad
	   g_{ab}^*(x)={\p^2 g(y)\over \p y^b \p y^a}
	   \big\vert_{y^a=S_1^a(x)}\,.
   \end{equation}

 Thus collecting the answers in  equations 
 \eqref{defofthickmorphism2a},\eqref{formalmapforinitialterm}
 and \eqref{formalmapforarbitraryterm} we come to
 \begin{proposition}\label{mapforthickmorphisms}
	 For thick morphism $\Phi_{S(x,q)}$
	 formal map $y^a(x)=y^a(x,g)$ in 
	 \eqref{defofthickmorphism2a} 
	 can be calculated 
	 recurrently by the equations
	 \eqref{formalmapforinitialterm},
	 \eqref{formalmapforarbitraryterm}.
	 In particular
up to  order  $k\leq 2$ it is defined by the following expression:
	      for arbitrary $g\in C^\infty(N)$,
	 $$
	 y^a(x)=y^a(x,g)=    
	 \underbrace {S_1^a(x)}_
	  {
	\hbox
	 {\footnotesize term of 
	 order $0$ in  $g$}}+
	 $$
      \begin{equation}\label{answerforformalmap1}
	      =
	      +\underbrace{2S_2^{ab}(x)g_b^*(x)}
	      _
	      {\hbox
	 {\footnotesize term of 
	 order $1$ in $g$}}+
	      \underbrace
	      {
3S^{abc}(x)g^*_b(x)g^*_c(x)+
	   4S^{ab}(x)S^{cd}(x)g^*_{bc}(x)g_d^*(x)
	      }_ {\hbox
	 {\footnotesize term of 
	 order $2$ in $g$}}
	      \, \,({\rm mod \A_3})\,.
      \end{equation}
 \end{proposition}	 
 Use this Proposition to 
   calculate  
   components $[\Phi_S^*(g)]_k$ of functional
   $\Phi_S^*(g)$.
    
   Due to definition \eqref{defofthickmorphism1} we have that 
                 $$
\Phi^*_S(g)=        
	\left( 
g(y^a)+S(x,q)-y^aq_a 
   \right)
\big\vert_{y^a={\p S(x,q)\over \p q_a}\,,
	   q_a={\p g(y)\over \p y^a}
	            }
		    =
		 $$
		\begin{equation}\label
		{calculationinthirdorder1}
	\left( 
g(y^a)+S_0(x)-\sum_{k\geq 2} (k-1)
	S_k^{a_1\dots a_k}(x)
{\p g(y)\over \p y^{a_1}}
\dots
{\p g(y)\over \p y^{a_k}}
   \right)
\big\vert_{y^a=y^a_0(x)+y^1_1(x)+\dots}\,,
		\end{equation}		
where 	
$y^a=y^a_0(x)+y^1_1(x)+\dots$  is a formal map 
\eqref{defofthickmorphism2a}.	
Here we used the fact that according to equations 
	\eqref{action1},
	\eqref{defofthickmorphism1} and
	\eqref{defofthickmorphism2}
	      $$
	      S(x,q)-y^aq_a=S(x,q)-{\p S(x,q)\over \p q_a}q_a=
	      \sum_k S_k^{a_1\dots a_k}(x)q_{a_1}\dots q_{a_k}-
	      \sum_k kS_k^{a_1\dots a_k}(x)q_{a_1}\dots q_{a_k}=
	      $$
	      $$
	      \sum_k (1-k)S_k^{a_1\dots a_k}(x)q_{a_1}\dots q_{a_k}\,.
	      $$
	  Now using equation \eqref{calculationinthirdorder1}
	 and equation 
  \eqref{answerforformalmap1} 
  in Proposition \ref{mapforthickmorphisms}
	  write down first few components $\left[\Phi_S^*(g)\right]_k$
of non-linear functional $\Phi_S^*(g)$      
    $$
    \left[  \left( g(y)\right)_{y^a(x)}\right]_{\leq 3}=
    g\left(y^a_{\leq 2}(x)\right)=
    g\left(y_0(x)+y^a_1(x)+y^a_2(x)\right)=
    $$
    $$
    g\left(S_1^a(x)+2S^{ab}(x)g_b^*(x)+
   4S^{ab}(x)S^{cd}(x)g^*_{bc}(x)g_d^*(x)
    \right)=
    $$
    $$
    g^*(x)+2S^{ab}(x)g_a^*(x)g_b^*(x)+
    3S^{abc}(x)g_c^*(x)g_b^*(x)g_a^*(x)+
    2S^{ab}(x)S^{cd}(x)g_{ab}^*(x)g^*_a(x)g^*_d(x)\,,
      $$
where we denoted by 
$\left[  g\left(y\right)_{y^a(x)}\right]_{\leq 3}$
 projection of functional 
 $g\mapsto g\left(y^a\left(x\right)\right)$ on $\A_3$.
 Hence it follows from  equation \eqref{calculationinthirdorder1}
 that
	 $$
\left[\Phi_S^*(g)\right]_{\leq 3}=
\left[\Phi_S^*(g)\right]_0+
\left[\Phi_S^*(g)\right]_1+
\left[\Phi_S^*(g)\right]_2+
\left[\Phi_S^*(g)\right]_3=
    $$
    $$
    S_0(x)+g\left(y^a_{\leq 2}(x)\right)
	-S_2^{ab}(x)
{\p g(y)\over \p y^{a}}
{\p g(y)\over \p y^{b}}
\big\vert_{y^a=S^a_0(x)+2S^{ab}(x)g_b^*(x)}-
		    $$
		    $$
-2S_3^{abc}(x)
{\p g(y)\over \p y^{c}}
{\p g(y)\over \p y^{b}}
{\p g(y)\over \p y^{a}}
\big\vert_{y^a=S^a_0(x)}		    
$$

 Collecting together the terms 
we come to formal power sums
we come to
\begin{proposition}
	\label{expressionforthickmorphism}
	Formal functional $\Phi_S^*(g)$	
	corresponding to thick morphism $\Phi_{S(x,q)}$
	can be calculated recurrently by equations
	\eqref{calculationinthirdorder1}.

In particular up to the order $\leq 3$ it is defined 
       	                 by the following expression 
			  $$
			\Phi_S^*(g)=
		       \underbrace{S_0(x)}_
		       {\footnotesize \hbox{term of order $0$ in $g$}}+
		     \underbrace
		           {
		       g\left(S^a(x)\right)
		       }_
		  {\footnotesize \hbox{term of order $1$ in $g$}}+
		       \underbrace
		       {
			       S^{ab}(x)g_b^*(x)g_b^*(x)
		       }_
		  {\footnotesize \hbox{terms of order $2$ in $g$}}+
                       $$
\begin{equation}\label{expansion2}
\underbrace
                            {
				    S^{abc}(x)
				    g_c^*(x)
				    g_b^*(x)
				    g_a^*(x)
		             +
			 2 S^{ac}
			  S^{bd}(x)g_{ab}^*(x)
			  g_d^*(x)
			  g_c^*(x)
	                       }_
{{\footnotesize \hbox{terms of order $3$ in $g$}}}
\, ({\rm mod} \A_4)\,.
\end{equation}
\end{proposition}

Thick morphisms define in general non-linear functionals 
$\Phi_S^*(g)$ belonging to space of formal functionals $\A$ (see
definition of formal functionals in 
\ref{defofformalfunctionals}).
As it was mentioned in   introduction
these non-linear functionals are non-linear homomorphisms.
Return to definition \ref{defofnonlinearhom1} of non-linear 
homomorphisms
formulating  it for formal functionals.
\begin{definition}\label{defoflocal}
	Let $L=L(x,g)$ be formal functional in $\A$
	(see definition \ref{defofformalfunctionals}). 
	According to definition \ref{defofnonlinearhom1}  
	this formal functional is {\it non-linear homomorphism}
	if its differential is usual
	homomorphism, i.e. for every function $g$ there exists 
a map 
     \begin{equation}\label{formalmap}
	y^a(x)=K^a(x,g)\,,
    \end{equation}
 such that for an arbitrary function $h$
	    \begin{equation}\label{differentialislocal}
	    L(g+\vare h)-L(g)=
	    \vare h \left(y^a(x,g)\right)\,,\quad(\vare^2=0)\,.
		    \end{equation}

	The  map $y^a(x,g)=K^a(x,g)$ 
in \eqref{formalmap}
	is in general a formal map:
       $$
	 y^a(x,g)=K^a_0(x)+
	 K^a_1(x,g)+   
	 K^a_2(x,g)+\dots=   
           $$
    \begin{equation}\label{formalmap1}
     K^a_0(x)+
     \int K^a_1(x,y)g(y)dy+
     \int K^a_1(x,y_1,y_2)g(y_1)g(y_2)dy_1dy_2+\dots
      \end{equation}
\end{definition}

Now we formulate

  \begin{theorem}\label{1}
Let $\Phi=\Phi_S\colon M\Rrightarrow  N$ be an arbitrary thick morphism.
	  Then formal functional  $\Phi_S^*(g)$
	  is non-linear homomorphism, i.e.
for arbitrary  functions $g$ 
	  there exists a map
	$y^a(x)=y^a(x,g)$    
	     such that for an arbitrary function
 $h$, ( $h\in C^\infty {N}$)
	  \begin{equation}\label{theorem1}
	    \Phi_S^*(g+\vare h)-
	    \Phi_S^*(g)=\vare h\left(y^a(x,g)\right)\,,\quad
		  \vare^2=0\,.
    \end{equation}
	 \end{theorem}
This very important observation was made by Voronov
in his pioneer work \cite {thickfirst} on 
thick morphisms. 

\begin{example}\label{newexample}
	For example consider pull-back
	   \begin{equation}\label{newexample1}
	  L(g)= \Phi_S^*(g)\,.
	   \end{equation}
	  According to Theorem \ref{1} this is
	  non-linear homomorphism. One can show that the map
	     $y^a=y^a(x,g)$ 
	     in equation \eqref{defofthickmorphism2a}
	     which we constructed above (see equations 
\eqref{formalmapforinitialterm}, \eqref{formalmapforarbitraryterm} 
and equation \eqref{answerforformalmap1} 
in Proposition \ref{mapforthickmorphisms})
	     is just formal map $K^a(x,g)$ \eqref{formalmap}
	     for this functional. 
	     (See the proof of Theorem \ref{1}
	     in the next section.)
\end{example}

For non-linear homomorphisms we will use
the notion of so called {\it support map}.
\begin{definition}\label{defofsupport}
If $L(g)$ is a functional which is non-linear homomorphism  then
a map $K_0^a(x)$ corresponding to the functional $L(g)$,
	which is the zeroth part of the formal map $K^a(x)$
	(see equations \eqref{formalmap} and \eqref{formalmap1}) 
will be called {\tt support  map corresponding to
functional $L(g)$}.

\end{definition}

\begin{example}
Consider functional	
$L(x,g)$ corresponding to thick morphism
	(see equation \eqref{newexample1}
	  in example \ref{newexample}).
If 
   $
	S(x,q)=S_0(x)+S_1^a(x)q_a+\dots	
	$
is generating function  \eqref{action1} which defines this 
	thick morphism,
	then it follows from equations \eqref{defofthickmorphism2a}
	and \eqref{formalmapforinitialterm} 
	that support map is equal to
	 $
	 K_0^a(x)=S^a_1(x)
	 $ (see also equation \eqref{answerforformalmap1}.)	 
\end{example}


\begin{definition}\label{associate}
	Let $L$ be an arbitrary functional in $\A$, 
	 \begin{equation*}
		 L(x,g)=\sum_k L_k(x,g)\,, \, {\rm where}\,\,
		 L_k(x,g)=[L(x,g)]_k\in \A_k
	 \end{equation*}
	(see equations \eqref{formalfunctional1} and 
	\eqref{formalfunctional1a}).  
	Taking the values of this functional on linear
	functions $y=y^al_a$ we assign to
	this functional, 
	{\it formal  function}
	 \begin{equation}\label{actionassociated1}
		 S_L(x,q)=L(x,g)\big\vert_{g=y^aq_a}=	 
	  S_0(x)+
			  \sum_k 
			  S_k^{a_1\dots a_k}(x)
			  q_{a_1}\dots q_{a_k}\,, \end{equation}
	where  
	tensors $\{S_k^{a_1\dots a_k}(x)\}$
	can be expressed through  polarised form of functionals
	 $L_k$ (see
	equations
	\eqref{polarisation1} and \eqref{polarisation2}
	in Appendix B):
	        \begin{equation*}
			S_k^{a_1\dots a_k}(x)=
			L_k^{\rm polaris.}
			\left(
			x,y^{a_1},\dots, y^{a_k}\right)\,,
	        \end{equation*}
	where $\{y^a\}$ are coordinates on $N$. E.g.
	            $$
	S_L^{ab}(x)=L_2^{\rm polaris.}(x,y^a,y^b)=
	 {1\over 2}\left(L_2\left(y^a+y^b\right)-
	 L_2\left(y^a\right)-L_2\left(y^a\right)\right)\,.
                $$
	We say that  $S_L(x,q)$ is 
	{\it formal function associated with
	functional $L$}.
	
\end{definition}

\smallskip

Let  $S=S(x,q)$ be an arbitrary formal generating function \eqref{action1}.
Let $\Phi_S$ be a thick morphism defined by 
this generating function, and let
$L(x,g)$ be a formal functional, $L(x,g)\in \A$, 
which defines pull-back of functions produced by
this thick morphism: $L(x,g)=\Phi^*_{S(x,q)}(g)$.
Then one can see that formal generating function
associated with functional $ L(x,g)=\Phi^*_{S(x,q)}(g)$
coincides with formal generating function $S(x,q)$:
	       \begin{equation}\label{actionisthesame1}
L(x,g)=\Phi^*_{S(x,q)}(g)\Rightarrow 	  
S_L(x,q)\equiv S(x,q)\,.     
	       \end{equation}
Indeed  in the case if function $g=y^al_a$ is linear then
calculations of pull-back $\Phi_S^*(g)$ 
by formulae \eqref{defofthickmorphism1}
and \eqref{defofthickmorphism2} become evident.
Indeed in this case we immediately come to 
equation \eqref{actionisthesame1} since
according
to equations \eqref{defofthickmorphism1}
and \eqref{defofthickmorphism2} 
    $$
 f(x)=g(y)+S(x,q)-y^aq_a=S(x,l)   
    $$
   because for linear function  $g(y)=y^aq_a$.
   
   It turns out that converse  implication is also valid
   for non-linear homomorphisms.

\begin{theorem}\label{2}
	Let $L=L(x,g)\in\A$ be an arbitrary
non-linear homomorphism, and let
	$S(x,q)$ be an action associated to it.  Then 
  \begin{equation*}
	  L(g)=\Phi_S^*(g)\,.
  \end{equation*}
\end{theorem}

   This is main result of this paper.

\section {Proof of the Theorems}
We recall here the proof of Theorem \ref{1}
and give a proof of Theorem \ref{2}.

\subsection {Proof of Theorem \ref{1}}
Check straightforwardly that 
a formal map $y^a(x,g)$ constructed 
in Proposition \ref{mapforthickmorphisms}  
(see equations 
\eqref{formalmapforinitialterm}, \eqref{formalmapforarbitraryterm} 
and equation \eqref{answerforformalmap1} 
in Proposition \ref{mapforthickmorphisms})
is just a map corresponding to function $g$ 
  i.e.	equation 
          \begin{equation}\label{checkstraightforwardly}
		  \Phi_S^*(g+\vare h)-  
		  \Phi_S^*(g)=\vare h(y(x,g))\,, (\vare^2=0)  
	  \end{equation} 
is satisfied.
(See also example \ref{newexample}.)

Using definition \eqref{defofthickmorphism1}
we see that in \eqref{checkstraightforwardly}
     $$
\Phi_S^*(g+\vare h)-  
\Phi_S^*(g)=
  $$
  $$
  \left[
	\left(g(y)+\vare h(y)\right)\big\vert_{y^a=y^a(x,g+\vare h)}+
S(x,q)\big\vert_{q_a=q_a(x,g+\vare h)}-
y^aq_a\big\vert_{y^a=y^a(x,g+\vare h),q_a=q_a(x,g+\vare h)}
\right]-
     $$
       \begin{equation}\label{checkstraightforwardly2}
  \left[
	  \left(g(y)\right)\big\vert_{y^a=y^a(x,g)}+
	       S(x,q)\big\vert_{q_a=q_a(x,g)}-
	       y^aq_a\big\vert_{y^a=y^a(x,g),q_a=q_a(x,g)}
\right]
     \end{equation}
     Here we introduced notation
       \begin{equation*}
	       q_a(x,g)={\p g(y)\over \p y^a}\big\vert_{y^a=y^a(x,g)}\,.
       \end{equation*}
      To see that right hand sides of equations
      \eqref{checkstraightforwardly}
      and \eqref{checkstraightforwardly2} coincide we note that
      in equation
      \eqref{checkstraightforwardly2} 
the following relations hold
\begin{equation}\label{A}
       	\left(g(y)+\vare h(y)\right)\big\vert_{y^a=y^a(x,g+\vare h)}
       	-
	g(y)\big\vert_{y^a=y^a(x,g)}=
	\vare{\p g(y)\over \p y^a}\big\vert_{y^a=y^a(x,g)}t^a=
	\vare q_a(x,g)t^a\,,
	 \end{equation}
	 \begin{equation}\label{B}
S(x,q)\big\vert_{q_a=q_a(x,g+\vare h)}-
S(x,q)\big\vert_{q_a=q_a(x,g)}=
\vare {\p S(x,q)\over \p q_a}\big\vert_{y^a=y^a(x,g)}r_a=
\vare y^a(x,g)r_a(x,g;h)\,,
	 \end{equation}
and
 \begin{equation}\label{C}
y^aq_a\big\vert_{y^a=y^a(x,g+\vare h),q_a=q_a(x,g)+\vare h)}-
y^aq_a\big\vert_{y^a=y^a(x,g,q_a=q_a(x,g)}=
\vare t^a q_a(x,g)+\vare y^a(x,g)r_a
\end{equation}
In equations \eqref{A}, \eqref{B} and \eqref{C}  we used notations
            $t^a,r_b$ such that
	        $$
y^a(x,g+\vare h)-		
y^a(x,g)=\vare t^a\,
{\rm and}\,\,		
q_a(x,g+\vare h)-		
q_a(x,g)=\vare r_a\,.
		$$
Comparing right hand sides of equations 
\eqref{A},
\eqref{B} and
\eqref{C} we come to conclusion that equation 
\eqref{checkstraightforwardly} is obeyed.\finish

\subsection {Proof of Theorem \ref{2}}
To prove Theorem \ref{2} we will formulate two lemmas.

\begin{lemma}\label{lemma1}

	Let $L=L(x,g)=\sum_{k\geq 0} L_k(x,g)$ 
	be an arbitrary functional
in $\A$
	which is non-linear homomorphism (see definition\ref{defoflocal}).
	Let $S_0(x)$ be a function which is equal to value
	of this functional on function $g=0$
	         \begin{equation}\label{affinecomponent}
	           S_0(x)=L(x,g)\big\vert_{g=0}\,,
		 \end{equation}	 
we will call sometimes this function 
	{\it an affine component of functional $L$}.
	
	Let a map
 $K_0^a(x)$ be a support  map corresponding to this functional
	(see definition \ref{defoflocal}\,).
Then
   \begin{equation*}
  L(g)=S_0(x)+g(K^a_0(x))\,\, ({\rm mod} \A_2)\ 
   \end{equation*}

\end{lemma}

\begin{lemma}\label{lemma2}
	Let $L(x,g)$ and $\widetilde L(x,g)$
	be two functionals on $\A$
	which both are non-linear homomorphisms,
and which coincide
up to the order $k-1$ ($k\geq 2$):
       \begin{equation*}
	       \begin{matrix}
		       \widetilde L(g)=		       
		       \sum_i\widetilde L_i(x,g)\,, 	       
		       \quad \widetilde L_i(x,g)\in A_i\cr 	       
		    L(g)=		       
		  \sum_i L_i(x,g)\,, 	       
		 \quad  L_i(x,g)\in A_i\cr 	       
		\widetilde L_j=L_j \,\, {\rm for}\,\,j\leq k-1\cr
	       \end{matrix}
       \end{equation*}
	Then 
the difference
of these functionals in the order $k$ is given by
$k$-linear functional $T_k(x,\p g)\in A_k$:
	\begin{equation*}
\widetilde  L_k(x,g)-L_k(x,g) =T_k(\p g)    
	\end{equation*}
where
           \begin{equation}\label{lemma2func}
\A_k\ni  T_k(\p g)=T^{a_1\dots a_k}(x)
		      g^*_{a_1}(x)\dots g^*_{a_k}
		      \quad {\rm and}\,\, 
   g_a^*(x)={\p g(y)\over \p y}\big\vert_
   {y^a=K^a(x)}\,,
	      \end{equation}
$K_0^a(x)$ is a support  map \ref{defofsupport} 
which is the same for both these functionals,
   and tensor $T^{a_1\dots a_k}$ is defined by  equation
			\begin{equation}
		\label{formulaforaction1}
T^{a_1\dots a_k}(x)=
\widetilde L_k^{\rm polaris.}
\left(x,y^{a_1},\dots, y^{a_k}
	\right)
-L^{\rm polaris.}_k
\left(x,y^{a_1},\dots, y^{a_k}\right)
			\end{equation}
			where  
$\widetilde L_k(x,g)$ and 
$L_k(x,g)$ 
are the terms of order $k$
in the expansion \eqref{formalfunctional1} 
of functionals $\widetilde L(x,g)$
and $L(x,g)$), and respectively
$\widetilde L^{\rm polaris.}_k(x,g_1,\dots,g_k)$ is
polarised form of functional $\widetilde L_k(x,g)$,
			and 
$L^{\rm polaris.}_k(x,g_1,\dots,g_k)$ is
polarised form of functional $L_k(x,g)$ 
(see equation \eqref{polarisation} 
in definition \ref{defofpolar} in Appendix B).

\end{lemma}

 Prove Theorem \ref{2} using these lemmas.

Let $L=L(g)$
be a functional in $\A$ which 
is non-linear homomorphism, i.e,
condition \eqref{differentialislocal} 
(see definition \ref{defoflocal}) holds 
for this functional, and  
     $$
  L(x,g)=L_0(x,g)+L_1(x,g)+\dots+L_k(x,g)+\dots\,,     
     $$
where every functional $L_r(x,g)$ has order $r$ in $g$: $L_r\in A_r$.

Consider an action $S(x,q)$ associated with this functional
(see equation \eqref{actionassociated1} in definition \ref{associate}).

 Consider the sequence of thick morphisms 
$\{\Phi_k\}$ ($k=0,1,2,\dots$)
  such that the thick morphism $\Phi_k$ is generated by the action
                        $$
     \bS_k(x,q)=
     S_0(x)+S_1^a(x)q_a+S_2^{ab}(x)q_aq_b+\dots+
S_k^{a_1\dots a_k}(x)q_{a_1}\dots q_{a_k}\,,
                        $$
and respectively the sequence $\{\Phi_k^*(g)=\Phi^*_{\bS_k}(g)\}$ 
of functionals, generated
by these thick morphisms.

Prove that for every $k$,
non-linear homomorphism
 $L(g)$ coincides up to terms of order $k$ in $g$ 
 with   functional $\Phi^*_{\bS_k}$:
             \begin{equation}\label{inductivestatement}
		    L(g)=\Phi_k^*(g) 
		    ({\rm mod }\, \A_{k+1})\,.
	      \end{equation}
This will be the proof of Theorem \ref{2}.

\begin{remark}   Thick morphisms $\{\Phi_k\}$ can be viewed
	as a sequense of morphisms tending to morphims $\Phi_S$.
\end{remark}

We prove equation \eqref{inductivestatement} by induction.
If $k=1$ then ${\bf S}_1(x)=S_0(x)+S_1^a(x)q_a$ and 
      $$
      \Phi_1^*(g)=
      S_0(x)+g(S_1^a(x))=L(g)
({\rm mod }\, \A_{2})\,.
      $$
due to Lemma \ref{lemma1}.  Thus equation \eqref{inductivestatement}
is obeyed if $k=1$.
Now suppose
that equation \eqref{inductivestatement} is obeyed for $k=m$, $m\geq 1$. 
Prove it for $k=m+1$.   Denote by
          \begin{equation}\label{temporarynotation}
\widetilde L(g)=\Phi_m^*(g)\,. 	  
	  \end{equation}  
Due to Theorem \ref{1} this functional is also 
non-linear homomorphism.
Both functionals are non-linear homomorphisms
and by inductive hypothesis functionals $L(g)$ and $\tilde L(g)$
coincide up to the order $m$. Hence lemma \ref{lemma2} implies 
that there exists tensor $T^{a_1\dots a_{m+1}}(x)$ such that 
	  \begin{equation}\label{comparing}
L(g)=\widetilde L(g)+T_{m+1}(\p g)= 
\Phi^*_{\bS_m}(g)+T_{m+1}(\p g) 
	  \,({\rm mod}\A_{m+2}),
	  \end{equation}
where 
        $$
T_{m+1}(\p g)=T^{a_1\dots a_{m+1}}(x)g_{a_1}^*(x)\dots g_{a_{m+1}}^*(x)\,,
\left(g_a^*(x)={\p g(y)\over \p y^a}\big\vert_{y^a=S^a_1(x)}\right)\,,
	  $$
and tensor $T^{a_1\dots a_{m+1}}(x)$ according
to equation \eqref{formulaforaction1} is defined by equation
	 \begin{equation}\label{defofT}
T^{a_1,\dots,a_{m+1}}_{m+1}=
L_{m+1}^{\rm polaris.}\left(x,y^{a_1},\dots,y^{a_{m+1}}\right)-
\widetilde 
L^{\rm polaris.}_{m+1}\left(x,y^{a_1},
\dots,y^{a_{m+1}}\right)\,,
	 \end{equation}
where $L^{\rm polaris.}_{m+1}$ is polarised form  of
functional $L_{m+1}(g)$ which contains terms of order
$m+1$ of functional $L(g)$.  Respectively functional
$\widetilde L^{\rm polaris.}_{m+1}$ is polarised form  of
functional $\widetilde L_{m+1}(g)$ which contains terms of order
$m+1$ of functional $\widetilde L(g)=\Phi_{\bS_m}(g)$.
It is easy to see that functional 
$\widetilde L^{\rm polaris.}_{m+1}$ is vanished
on arbitrary linear functions:
	  \begin{equation}\label{vanishesonlinear}
		  \widetilde L\left(x,l_1,\dots,l_{m+1}\right)=0\,,\quad
		  \hbox{if functions $l_i$ are linear: 
		  $l_i=y^al_{ai}$, $i=1,\dots,m+1$}\,.§:
          \end{equation}
Indeed functional $\widetilde L(g)=\Phi^*_{\bS_m}(g)$
is assigned to the action $\bS_m(x,q)$
 which is a polynomial of order $\leq m$, hence
due to equation \eqref{actionisthesame1}
it vanishes for arbitrary linear function $g=y^al_a$,
hence polarised form vanishes also on linear functions
( see equation \eqref{polarisation1} in Appendix B).
Thus we come to condition  \eqref{vanishesonlinear}.
This condition means that in particular 
      $$
\widetilde 
L^{\rm polaris.}_{m+1}\left(x,y^{a_1},
\dots,y^{a_{m+1}}\right)=0\,,\quad{\rm for}\,\, 
\widetilde L(g)=\Phi_{m+1}^*(g)\,,
      $$
      hence we come to conclusion that
 tensor
  $T^{a_1\dots a_{m+1}}(x)$ in equation \eqref{defofT}
  is equal to $S^{a_1\dots a_{m+1}}(x)$.

 We see that 
   \begin{equation}\label{m+2before}
L(g)=\Phi^*_m(g)+S_{m+1}(\p g)\,\, 
	({\rm mod} \A_{m+2})\,.
	      \end{equation}
	    On the other hand up to the terms of
	    order $m+1$, right hand sight of this equation
	   is equal to $\Phi^*_{m+1}$: 
   \begin{equation}\label{m+2after}
\Phi^*_{m+1}(g)=\Phi^*_m(g)+S_{m+1}(\p g) 
	({\rm mod} \A_{m+2})\,.
	      \end{equation}
One can see it straightforwardly using equation
\eqref{defofthickmorphism1} or  it is much easier to check 
equation taking
differential of this equation.  Namely taking differential
of equation \eqref{m+2after} and using equations 
\eqref{formalmapforarbitraryterm}
and \eqref{checkstraightforwardly} we come to equation
                $$
h\left(y^a_{m+1}(x,g)\right)		
=h\left(y^a_{m}(x,g)\right)+S_{m+1}^{aa_1\dots a_m}
    g_{a_1}^*\dots
    g_{a_m}^*
	({\rm mod} \A_{m+1})\,,
		$$
		where 
$y^a_{\bS_{k}}(x,g$ is a map $y^a(x,g)$ corresponding
to thick morphism $\Phi_{\bS_k}$
($\Phi^*_{\bS_k}(g+\vare h)-\Phi^*_{\bS_k}(g)=
h\left(y^a_{\bS_{m+1}}(x,g)\right) h$).
	Comparing left hand sides of equations
	\eqref{m+2before} and
	\eqref{m+2after} we see that
equation \eqref{inductivestatement} 
	holds for $k=m+1$.
	    This  ends the proof.

	    \finish.   

It remians to prove lemmas.

\section {Proofs of lemmas}

\subsection  {Proof of the Lemma \ref{lemma1}}

Let $L=L(x,g)$ be a functional in $\A$  which is non-linear homomorphism.
      \begin{equation}\label{lemma11}
   L(x,g)=
	      L_0(x)+L_1(x,g)+\dots=
	      L_0(x)+L_1(x,g)({\rm mod\,} \A_2)   
      \end{equation}
If we put $g=0$ we come to 
    $L_0(x)=S_0(x)=L(g)\big\vert_{g=0}$.

Differentiate equation \eqref{lemma11}. 
Using equation \eqref{formalmap1} we come to
          $$
	  L(x,g+\vare h)-L(x,g)=
	  \vare 
	  h\left(
	  y^a(x,g)
	  \right)=
	  \vare h
	 \left(
	  K_0^a(x)+K_1^a(x,g)+\dots
	  \right)
	   =\vare h\left(
	  K_0^a(x)
	  \right) \, ({\rm mod} \A_1)
    $$
This is true for arbitrary smooth function $h$. This implies that
  $L_1(x,g)=g\left(K_0^a(x)\right) $. Hence 
         $$
 L(g)=L_0(g)+L_1(g) \, ({\rm mod} \A_2)=
 S_0(x)+g\left(K_0^a(x)\right)\,({\rm mod} \A_2)=
	 $$
	 $$
S_0(x)+\int K(x,y)g(y)dy+\hbox{\footnotesize terms of order 
$\geq 2$ in $g$}\,,\quad {\rm with}
   \,\, K(x,y)=\delta(y^a-K^a_0(x))\,.
          $$

First lemma is proved.

\subsection {Proof of lemma \ref{lemma2}}

Let functionals $L(g)$ and $\widetilde L(g)$
both be functionals which are non-linear homomorphisms
(see definition \ref{defoflocal}).
Suppose these functionals  coincide up to the order
$k-1$ ($k=2,3,\dots$). According to expansion
\eqref{formalfunctional1} this means that difference of these functionals
is a functional $T_k(g)$ of order $k$
               \begin{equation}\label{lemma21}
		\widetilde L(g)-L(g)=T_k(x,g)\in \A_{k+1}\, \quad
		       {\rm i.e. }\,\,
	\widetilde L(g)-L(g)-T_k(g)=0	({\rm mod} \A_{k+1})\,,
	\end{equation}
where 
	   \begin{equation}\label{lemma22}
		   T_k(x,g)=\int T(x,y_1,\dots,y_k)
		   g(y_1)\dots g(y_k)dy_1\dots dy_k\,.   
	   \end{equation}
Take the differential of  equation \eqref{lemma21}. We come to
$$
\left(\widetilde L(g+\vare h)-	  
\widetilde L(g)\right)-
   \left(
 L(g+\vare h)-	  
 L(g)\right)=
 \vare h\left(\widetilde {y^a}(x,g)\right)	  
		  -\vare h\left( {y^a}(x,g)\right)= 
	  $$
 $$
=T_k(x,g+\vare h)-T_k(g)+\hbox {terms of order $\geq k$ in $g$}=
   $$ 
	  \begin{equation}\label{lemma23}=
		  \vare kT^{\rm polaris.}_k
 \left(h,\underbrace {g,\dots,g}_{\hbox {$k-1$ times}}\right)+
       \hbox {terms of order $\geq k$ in $g$}=
          \end{equation}
Here $T^{\rm polaris.}_k=T_k(x,g_1,g_2,\dots,g_k)$ is the
polarisation of the form $T_k(x,g)$ (see equation \eqref{polarisation1}
in definition \ref{defofpolar}).    Recall that  
if function $T(x,y_1,y_2,\dots,y_k)$
which correspond to functional $T_k(g)$ in equation \eqref{lemma22}
is symmetric function on variables $y_1,\dots y_k$ then
                    (see equation \eqref{generalisedfunctions})  
                    $$
       T^{\rm polaris.}_k(x,g_1,\dots,g_k)=
		  \int T(x,y_1,y_2,\dots,y_k)
		   g_1(y_1)g_2(y_2)\dots g_k(y_k)dy_1 dy_2\dots dy_k\,.
	                $$

Formal maps 
$y^a(x,g)$  corresponding to differential
   $dL(g)=L(g+\vare h)-L(g)$ of functional $L(g)$
   and 
$\widetilde y^a(x,g)$  corresponding to differential
   $d\widetilde L(g)=
\widetilde L(g+\vare h)-\widetilde L(g)$ of functional
   $\widetilde L(g)$ 
according to equation \eqref{formalmap1}
are given by formal power series
      $$
   y^a(x,g)=
   K^a_0(x)+   
   K^a_1(x,g)+\dots+   
     K^a_{k-2}(x,g)+K^a_{k-1}(x,g)+
{\footnotesize \hbox {terms of order $\geq k$ in $g$}}
     $$
and
      \begin{equation}\label{lemma25}
   \widetilde y^a(x,g)=
   \widetilde K^a_0(x)+   
   \widetilde K^a_1(x,g)+\dots+   
    \widetilde  K^a_{k-2}(x,g)+
    \widetilde K^a_{k-1}(x,g)+ 
{\footnotesize \hbox {terms of order $\geq k$ in $g$}}\,.
     \end{equation}
  Recall that    here $K^a_r(x,g)$ and $\widetilde K^a_r(x,g)$
are maps of order $r$ in $g$:
            $$
K^a_r(x,g)=\int K(x,y_1,\dots, y_r)g(y_1)\dots g(y_r)
dy_1\dots dy_r\,.
      $$
Since  functionals $L(g)$ and $\widetilde L(g)$ 
coincide up to the order $k-1$, their differentials
coincide up to the order $k-2$. 
Hence it follows from equation \eqref{lemma23}
that in equation \eqref{lemma25} all the maps $K^a_r$ coincide with 
maps $\widetilde K^a_r$ for $r=0,1,2,\dots, k-2$
           \begin{equation*}
K_0^a(x)	   
=\widetilde K_0^a(x)\,,
\dots\,,	   
K_{k-2}^a(x,g)	   
=\widetilde K_{k-2}^a(x,g)\,,
	   \end{equation*}
and it is the 
difference between maps $\widetilde K_{k-1}$ and $K_{k-1}$ which
produces the functional $T_k(x,g)$. 

Rewrite equation 
 \eqref{lemma23} projecting all terms on subspace $A_{k-1}$.
 We come to
	 $$
\left[\left(\widetilde L(g+\vare h)-	  
\widetilde L(g)\right)\right]_{k-1}-
   \left[
 L(g+\vare h)-	  
 L(g)\right]_{k-1}=
 \vare 
    \left[h\left(\widetilde {y^a}(x,g)\right)	  
		  -\vare h\left( {y^a}(x,g)\right)
		  \right]_{k-1}= 
       $$
 $$
{\p h\over \p  y^a}\big\vert_{y^a=K^a_0(x)}
\left[ \widetilde K_{k-1}^a(x,g)-
   K_{k-1}^a(x,g)\right]
           =
{\p h\over \p  y^a}\big\vert_{y^a=K^a_0(x)}
P^a_{k-1}(x,g)=
  $$
  $$
  $$
 $$
=T_k(x,g+\vare h)-T_k(g)=
   $$ 
	  \begin{equation*}\label{lemma26}=
		  \vare kT^{\rm polaris.}_k
 \left(x,h,\underbrace {g,\dots,g}_{\hbox {$k-1$ times}}\right)
		  \,.
          \end{equation*}
where we denote by $P^a_{k-1}(x,g)$ the difference between maps
   $\widetilde K_{k-1}^a(x,g)$
  and  $ K_{k-1}^a(x,g)$
    \begin{equation*}\label{lemma26}
 P^a_{k-1}(x,g)=\widetilde K_{k-1}^a(x,g)-
   K_{k-1}^a(x,g)=\int P^a_{k-1}(x,y_1,\dots,y_{k-1})
      g(y_1)\dots g(y_{k-1})dy_1\dots dy_k\,.
   \end{equation*}
	The map $P^a_{k-1}(x,g)$ has order $n-1$ over $g$.
	Consider polarisation 
	$P^{a\, {\rm polaris.}}_{k-1}(x,_1,\dots,g_{k-1})$  
	\eqref{polarisation1} of  this map.
	Equation \eqref{lemma26} implies 
	\begin{equation*}
{\p h\over \p  y^a}\big\vert_{y^a=K^a_0(x)}
P^{a\, {\rm polaris.}}_{k-1}(x_1,\dots,g_{k-1})  
		  =  \vare kT^{\rm polaris.}_k
		  \left(x,h,g_1,\dots,g_{k-1}\right)
		\big\vert_{g_1=\dots=g_{k-1}=g}
		  \,.
	\end{equation*}
Thus we come to equation
\begin{equation}\label{lemma27}
T^{\rm polaris.}_k
 \left(x,g_1,\dots,g_{k}\right)=
 {1\over k}{\p g_1\over \p  y^a}\big\vert_{y^a=K^a_0(x)}
	P^{a\,{\rm polaris.}}_{k-1}(x,g_2,\dots,g_k)\,,
\end{equation}
where $g_1,\dots,g_k$ are arbitrary functions
and left hand side of this equation is symmetric with
 respect to transposition of functions $\{g_1,\dots,g_k\}$.
It follows from equation \eqref{lemma27} that
            $$
	P^{a\,{\rm polaris.}}_{k-1}(x,g_2,\dots,g_k)\,,
	 = 
kT^{\rm polaris.}_k
 \left(x,y^a,g_2,\dots,g_{k}\right)
	  $$
hence
\begin{equation}\label{lemma2key}
T^{\rm polaris.}_k
 \left(x,g_1,\dots,g_{k}\right)=
{\p g_1\over \p  y^a}\big\vert_{y^a=K^a_0(x)}
T^{\rm polaris.}_k
 \left(x,y^a,g_2,\dots,g_{k}\right)\,.
\end{equation}
   Equation \eqref{lemma2key} and symmetricity of 
 functional $T_k(x,g_1,\dots,g_k)$ imply that
    $$
T^{\rm polaris.}_k
 \left(x,g_1,g_2,\dots,g_{k}\right)=
 {\p g_1\over \p  y^a}\big\vert_{y^a=K^a_0(x)}
T^{\rm polaris.}_k
 \left(x,y^a,g_2,\dots,g_{k}\right)=
 $$
 $$
	    T^{\rm polaris.}_k
 \left(x,g_2,g_1,\dots,g_{k}\right)=
 {\p g_2\over \p  y^a}\big\vert_{y^a=K^a_0(x)}
T^{\rm polaris.}_k
 \left(x,y^a,g_1,\dots,g_{k}\right)=\dots=
$$
$$
{\p g_1\over \p  y^{a_1}}\big\vert_{y^{a_1}=K^{a_1}_0(x)}
\dots
{\p g_k\over \p  y^{a_k}}\big\vert_{y^{a_k}=K^{a_k}_0(x)}
 T^{\rm polaris.}_k
 \left(x,y^{a_1},\dots,y^{a_k}\right)=
$$
\begin{equation}\label{mainstatement}
	=g_{a_1}^*(x)	
	\dots g_{a_k}^*(x)
	T^{a_1\dots a_k}(x)\,,
\end{equation}
where 
    $$
   T^{a_1\dots a_k}(x)=T_k\left(x,y^{a_1},\dots,y^{a_k}\right)=
   {\rm and}\,\,
g_a^*(x)=
{\p g\over \p  y^{a}}\big\vert_{y^{a}=K^{a}_0(x)}\,.
$$
Now returning to equation \eqref{lemma21} and comparing it 
with formulation of lemma \ref{lemma2} we come to proof of 
lemma \ref{lemma2}:
      $$
    \widetilde L_k(g)-L_k(g)=T_k(x,g_1,\dots,g_k)
    \big\vert_{g_1=\dots=g_k=g}=T_k(\p g)\,. 
      $$

 \def\cprime{$'$}

\medskip

{\small
\section{ Appendix A. Thick morphisms and $L_\infty$ maps}


We briefly here discuss why thick morphisms 
is an adequate tool to describe $L_\infty$-morphisms 
of homotopy Poisson algebras
(see \cite{thickfirst} and  \cite{thicktheor} for detail).
  For this purpose we need to consider thick morphisms of supermanifolds.
  However we can catch some improtant features considering
  just usual manifolds.
  We first consider thick morphisms  
  for usual manifolds, and show that
  in this case thick morphisms describe morphisms  of algebras 
  of functions on these manifolds which are provided with multilinear
  symmetric brackets. It turns out that if we consider supermanifold, then
  under some assumptions these algebras become homotopy Poisson 
  algebras. 

	Let  $M$ be an arbitrary manifold, and $H=H(x,p)$
	be a function (Hamiltonian) on cotangent bundle $T^*M$.
	This Hamiltonian  $H$ defines the series of symmetric
	brackets on $M$ via canonical symplectic structure on $T^*M$ 
	          $$
		  \langle \emptyset \rangle_H\,,
		  \langle f_1 \rangle_H\,,
		  \langle f_1, f_2 \rangle_H\,,
		  \langle f_1, f_2,f_3 \rangle_H\,,
		  \dots
		  \langle f_1, f_2,\dots,f_k \rangle_H\,,
		  $$
  
where
         $$
 \langle \emptyset \rangle_H=H(x,p)\big\vert_{p=0}=H_0(x)
          $$
	  $$
 \langle f_1 \rangle_H=\left(H,f_1\right)\big\vert_{p=0}=
  H_1^a(x){\p f_1(x)\over \p x^a}\,,
	 $$
	  $$
 \langle f_1,f_2 \rangle_H=
 \left(\left(H,f_1\right),f_2\right)\big\vert_{p=0}=
  H_1^{ab}(x)
  {\p f_1(x)\over \p x^a}
  {\p f_2(x)\over \p x^a}\,,
	 $$
and so on:
	  \begin{equation}\label{homotopybracketsformanifoldsusual}
 \langle f_1,f_2,\dots,f_k \rangle_H=
 \underbrace
 {(\dots (}
 _{\hbox {$k$ times}}
 H,f_1),f_2 )\dots f_k)\big\vert_{p=0}=
  H_k^{a_1\dots a_k}(x)
  {\p f_1(x)\over \p x^{a_1}}
  \dots
  {\p f_k(x)\over \p x^{a_k}}\,.
        \end{equation}
Here   $(\_,\_)$ is Poisson bracket on $T^*M$ 
corresponding to canonical symplectic structure:
	     \begin{equation}\label{canonicalpoissonbracket}
	     \left(f(x,p),g(x,p)\right)=
	     {\p f(x,p)\over \p p_a}
	     {\p g(x,p)\over \p x^a}
	     -
	     {\p g(x,p)\over \p p_a}
	     {\p f(x,p)\over \p x^a}\,.
	     \end{equation} 
We suppose that 
Hamiltonian  $H=H(x,p)$ is a {\it formal Hamiltonian}, i.e.
formal function, power series over $p$:
	       \begin{equation*}\label{formalhamiltonian}
H=H(x,p)=H_0(x)+H_1^a(x)p_a+H_2^{ab}(x)p_bp_a+
H_3^{abc}(x)p_cp_bp_a+\dots
	       \end{equation*} 
where  all coefficients are smooth functions on $x$.

\begin{remark}\label{canonicalcoordinatesremark}
All  these formulae are written in local coordinates $(x^a,p_b)$
in $T^*M$ corresponding to local coordinates $x^a$ on $M$
(if $x^{a'}$ are new local coordinates on $M$,
then new local coordinates $(x^{a'}, p_{b'})$) are
		 \begin{equation}\label{canonicalcoordinatesonbundle}
    x^{a'}=x^{a'}(x)\,,  p_{b'}={\p x^b(x')\over \p x_{b'}}p_b\,.
		   \end{equation}
\end{remark}  

Notice that every  Hamiltonian $H(x,p)$
defines vector field 
	    \begin{equation*}\label{vectorfieldonfunctions1}
		    X_H=\int 
		   H\left(
		    f\left(x\right),
		    {\p f\left(x\right)\over \p x}
		     \right)dx    
	    \end{equation*} 
on the space of function. Vector field $X_H$  assigns to every function
 $f\in C^\infty(M)$  infinitesimal curve
  	    \begin{equation}\label{vectorfieldonfunctions2}
		   f+\vare X_H= 
	         f(x)+
		   \vare H\left(
		    f\left(x\right),
		    {\p f\left(x\right)\over \p x}\right)\,,
		    \quad (\vare^2=0)\,.
	    \end{equation} 
  \smallskip

  Now consider two  manifolds   $M$ and  $N$. Let $H_M(x,p)$
  be formal Hamiltonian on $M$,
  and let $H_N(y,q)$
  be formal Hamiltonian on $N$.  
  Hamiltonian  $H_M(x,p)$  induces on $M$
  the sequence of multilienar symmetric brackets
     $\left\{\langle f_1,\dots,f_p\rangle_M\right\}$ 
     on functions on $M$,  
     and respectively
   Hamiltonian  $H_N(y,q)$  induces on $N$
  the sequence of multilinear symmetric brackets
     $\left\{\langle g_1,\dots,g_q\rangle_M\right\}$ 
     on functions on $N$
     ($p,q=0,1,2,3,\dots$).

  We say that formal functional $L(g)$
  is morphism of multilinear symmetric brackets
  on $N$ to  multilinear symmetric brackets on $M$
  if vector fields $X_{H_M}$ and $X_{H_N}$
  are connected by functional $L(g)$, i.e.
  according to formulae \eqref{vectorfieldonfunctions2}
      \begin{equation*}\label{bracketsconnected}
       L\left(g+\vare X_N\right)=
       L\left(g\right)+\vare X_M\,.
      \end{equation*} 
   Consider thick morphism 
  $\Phi_S\colon M\Rrightarrow N$  generated
  by $S(x,q)$ 
  and consider  formal functional 
  $\Phi_S^*(g)$ on $C^\infty(N)$ defined by this thick morphism
(see equations \eqref{action1}---\eqref{expansion2}
and remark \ref{infactmanifolds}).

We say that Hamiltonians $H_M$ and $H_N$ are $S$-related
if 
   	    \begin{equation*}\label{hamiltoniansrelated}
		    H_M\left(x,{\p S\left(x,q\right)\over \p x}\right)
		    \equiv
		   H_N\left({\p S\left(x,q\right)\over \p q},q\right)
	    \end{equation*} 
	    The following remarkable theorem takes place:

\begin{theorem}\label{Voronovtoy}(Voronov, 2014)	    
If  Hamiltonians $H_M$ and $H_N$ are  $S$-related, then
	formal functional $L(g)$ defined by thick morphism
	$\Phi_S$,   $L(g)=\Phi_S^*(g)$
defines morphisms of multilinear brackets
 
     $\left\{\langle f_1,\dots,f_p\rangle_M\right\}$ 
	and
     $\left\{\langle g_1,\dots,g_q\rangle_M\right\}$ 
     $\left\{\langle g_1,\dots,g_q\rangle_M\right\}$.
In other words thick morphism connects
these brackets.
\end{theorem}

\smallskip

Now consider the case of supermanifolds.

In this case all the constructions above will remain the same,
	just in some formulae will appear a sign factor.
	(See \cite{thickfirst} and  \cite{thicktheor} for detail).
In particular arbitrary Hamiltonian $H=H(x,p)$ which is
a function on cotangent bundle $T^*M$ to supermanifold $M$
will define  the collection of symmetric brackets like in the
case \eqref{homotopybracketsformanifoldsusual}. On the other hand 
if Hamiltonian $H_M$ is {\it odd} 
and Hamiltonian $H_M$ obeys condition
	   \begin{equation}\label{master1}
		   \left(H_M,H_M\right)\equiv 0\,,
	   \end{equation}
then these brackets will become {\it  homotopy Poisson brackets}.
This is famous construction of homotopy Poisson brackets
derived by odd Hamiltonian $H_M$ which obeys so called
master-equation \eqref{master1} (see for detail \cite{thickkoszul}).

\medskip

\section{Appendix B. Polarisation of functionals}  
 
 It is useful 
 to consider polarised form of formal functionals.

  \begin{definition}\label{defofpolar}
	  Let $L_k(x,g)$ be formal functional
	  of order $k$, $L_k(x,g)\in \A_k$
	  (See for definition \ref{defofformalfunctionals}.)
	  Polarisation
	  of functional $L_k(x,g)$  
	  is 
	  the functional $L_k^{\rm polaris.}(x,g_1,\dots,g_k)$
	  which linearly depends on $k$ functions
	  $g_1,\dots,g_k$
       such that  for every function $g$
	     \begin{equation}\label{polarisation1}
	L_k(x,g) =
	L_k^{\rm polaris.}(x,g_1,\dots, g_k)
	\big\vert_{g_1=g_2=\dots=g_k=g}\,.   
	     \end{equation}
	 Using elementary combinatoric one can express
	  polarised form  
	$L_k^{\rm polaris.}(x,g_1,\dots, g_k)$
	  explicitly 
	  in terms of functional $L_k(x,g)$, ($L_k\in A_k$):
	  \begin{equation}\label{polarisation2}
 L_k^{\rm polaris.}(x,g_1,\dots,g_k)=
		  {1\over k!} \sum(-1)^{k-n}
		  L_k\left(
		  x,g_{i_1}+\dots+g_{i_n}	  
		  \right)\,,	 
	 \end{equation}
where summation goes over all  non-empty subsets
	  of the set $\{g_1,\dots,g_k\}$. E.g.
	   if $L=L_3$ then
	   $$
	L^{\rm polaris.}(x,g_1,g_2, g_r)=
	  {1\over 6}
	  \left(
	  L_3\left(x,g_1+g_2+g_3\right)
	  -L_3\left(x,g_1+g_2\right)
	  -L_3\left(x,g_1+g_3\right)
	  -L_3\left(x,g_2+g_3\right)\right.
	     $$
	     $$
	     \left.
	  +L_3\left(x,g_1\right)
	  +L_3\left(x,g_2\right)
	  +L_3\left(x,g_3\right)
	  \right)\,.
	    $$
    If functional $L_r(x,g)$ is 
	  expressed through (generalised) functions
	  $L(x,y_1,\dots,y_r)$	
(see equation \eqref{generalisedfunctions}) such that it  is symmetric 
	  with respect to coordinates 
	  $y_1,\dots,y_r$
	 then
	  \begin{equation}\label{polarisation}
 L^{\rm polaris.}(g_1,\dots,g_r)=
\int L(x,y_1,\dots,y_r)g_1(y_1)\dots g(y_r)dy_1\dots dy_r\,.	 
	 \end{equation}
\end{definition} 
It is useful also to note that if 
$L(x,g)=L_0(x)+L_1(x,g)+\dots+L_n(x,g)$ then for every
$k\colon\, k=0,1,\dots,n$
     	  \begin{equation}\label{polarisation3}
 L_k^{\rm polaris.}(x,g_1,\dots,g_k)=
		  {1\over k!} \sum(-1)^{k-n}
		  L\left(
		  x,g_{i_1}+\dots+g_{i_n}	  
		  \right)\,,	 
	 \end{equation}
where summation goes over all  subsets
	  of the set $\{g_1,\dots,g_k\}$
including empty subset. 
(For empty subset $L(x,\emptyset)=L_0(x)$.)	  
 }

\end{document}